\documentclass[11pt]{article}
\usepackage{amscd,amssymb}

\usepackage[bookmarksnumbered,linkcolor=black,citecolor=black]{hyperref}

\newtheorem{thm}{Theorem}[section]
\newtheorem{corol}[thm]{Corollary}
\newtheorem{lemma}[thm]{Lemma}

\newtheorem{defi}[thm]{Definition}
\newtheorem{remark}[thm]{Remark}
\newtheorem{example}[thm]{Example}
\newenvironment{proof}{\noindent {\sl Proof.}}{$\Box$ \bigskip}
\def\Z {{\Bbb Z}}

\title  {Minimality in topological groups and Heisenberg type groups}
\author{Menachem Shlossberg\\ \footnotesize{Department of Mathematics,
Bar-Ilan University, 52900 Ramat-Gan, Israel}\\
\footnotesize{email: {\tt menysh@yahoo.com}}}
\date{June 14, 2007}

\begin{document}

\maketitle

\begin{abstract}
We study relatively minimal subgroups in topological groups. We
find, in particular, some natural relatively minimal subgroups in
unipotent groups which are defined over "good" rings. By "good"
rings we mean archimedean absolute valued (not necessarily
associative) division rings. Some of the classical rings which we
consider besides the field of reals are the ring of quaternions and
the ring of octonions. This way we generalize in part a previous
result which was obtained by Dikranjan and Megrelishvili
\cite{DIM04} and involved the Heisenberg group.
\end{abstract}

\section{Introduction}
A Hausdorff topological group $G$ is minimal if $G$ does not admit a
strictly coarser Hausdorff group topology or equivalently   if every
injective continuous group homomorphism $G\rightarrow P$ into a
Hausdorff topological group is a topological embedding. The concept
of minimal topological groups was introduced by Stephenson
\cite{STE71}
and Do\"\i chinov \cite{DOI72} in 1971 as a natural generalization
of compact groups.

Heisenberg group  and more precisely its generalization, which we
present in section ~\ref{sec:mgh}  (see  also \cite{MIL89,REI74}),
provides many examples of minimal groups.

Recently Dikranjan and Megrelishvili \cite{DIM04} introduced the
concept of {\it co-minimality} (see Definition ~\ref{def:com}) of
subgroups in topological groups after the latter author had
introduced the concept of {\it relative minimality} (see Definition
~\ref{def:rel} and also \cite{MEG04})  of subgroups in topological
groups and found such subgroups in a {\it generalized Heisenberg
group} (see \cite{MIL89,REI74}).

In \cite[Proposition 2.4.2]{DIM04} Megrelishvili and Dikranjan
proved that the canonical bilinear mapping $V\times
V^{\ast}\rightarrow \Bbb R,\ < v,f >= f(v)$ is \emph{strongly
minimal} (see Definition ~\ref{def:stm}) for all 
normed spaces $V.$

The following result is obtained as a particular case:  The inner
product map $$\Bbb R^n \times \Bbb R^n \rightarrow \Bbb R$$  is
strongly minimal.
The latter result leads in \cite{DIM04} and \cite{MEG04} to the
conclusion that for every $n\in  \Bbb N $ the
subgroups $$\bigg\{\left (\begin{array}{ccc} 1 & a& 0\\
0 & I_n & 0\\
0 & 0 & 1\\
 \end{array}\right)\bigg |\ a\in \Bbb R^n \bigg \}, \ \bigg\{\left (\begin{array}{ccc} 1 & 0 & 0\\
0 & 1 & b\\
0 & 0 & 1\\
 \end{array}\right)\bigg |\ b \in \Bbb R^n \bigg \} $$ are relatively
minimal in the group  $$\bigg\{\left (\begin{array}{ccc} 1 & a& c\\
0 & I_n & b\\
0 & 0 & 1\\
 \end{array}\right)\bigg |\ a,b\in \Bbb R^n, \ c\in \Bbb R \bigg \}$$ 
  which is known as  the classical $2n+1$-dimensional
Heisenberg group (where $I_n$ denotes the identity matrix of size
$n$). Theorem ~\ref{thm:new} and Corollary ~\ref{cor:nrs} generalize
these results and allow us to replace the field of reals by every
other archimedean absolute valued (not necessarily associative)
division ring, for example, they can be applied for the ring of
quaternions and the ring of octonions.  Theorem ~\ref{thm:ne2}
provides a different generalization. It generalizes the case of the
classical real $3$-dimensional Heisenberg group. We consider for
every $n\in \Bbb N$ the group of upper unitriangular matrices over
an archimedean absolute valued field of size $n+2\times n+2$ and we
find relatively minimal subgroups of this group. This result is a
generalization since the classical real $3$-dimensional Heisenberg
group is a unitriangular group. This theorem is not new when we take
$n=1$ and consider the field to reals. However, we obtain a new
result even for $\Bbb R$ when we take $n>1.$ This theorem can also
be applied for the fields $\Bbb Q$ and $\Bbb C.$

 \section{Minimality in Generalized Heisenberg
groups}\label{sec:mgh} \vskip 0.3cm
The group $$H=\bigg\{\left (\begin{array}{ccc} 1 & x & a\\
0 & 1 & y\\
0 & 0 & 1\\
 \end{array}\right)\bigg |\  x,y,a \in \Bbb R \bigg \} \cong (\Bbb
 R \times \Bbb
 R)\leftthreetimes \Bbb  R$$
is known as the classical real 3-dimensional  {\it Heisenberg Group}.\\

\vskip 0.5cm We need a far reaching generalization \cite{MIL89,
REI74, MEG04}, the {\it generalized Heisenberg group}, which is
based on biadditive mappings.
\begin{defi}
Let $E,F,A$ be abelian groups. A map $w:E\times F \rightarrow A$ is
said to be {\it biadditive }if the induced mappings
$$w_{x}:F\rightarrow A,\ w_{f}:E\rightarrow A,
\ w_{x}(f):=w(x,f)=:w_{f}(x)$$ are homomorphisms for all $x\in E$
and $f\in F$.
\end{defi}
\begin{defi}\cite[Definition 1.1]{MEG04}
Let $E,F$ and $A$ be Hausdorff abelian topological groups and
$w:E\times F \rightarrow A$ be a continuous biadditive mapping.
Denote by $H(w)=(A\times E)\leftthreetimes F$ the topological
semidirect product (say, {\it generalized Heisenberg group} induced
by $w$) of $F$ and the group $A\times E$. The group operation is
defined as follows: for a pair
$$u_{1}=(a_{1},x_{1},f_{1}), \ u_{2}=(a_{2},x_{2},f_{2}) $$ we
define $$u_{1}u_{2}=(a_{1}
+a_{2}+f_{1}(x_{2}),x_{1}+x_{2},f_{1}+f_{2})$$ where,
$f_{1}(x_{2})=w(x_{2},f_{1})$. Then $H(w)$ becomes a Hausdorff
topological group. In the case of a normed space $X$ and a canonical
biadditive function $w:X\times X^{\ast}\rightarrow \Bbb R $
$(x,f)\mapsto f(x)$ (where $X^{\ast}$ is the Banach space of all
continuous functionals from $X$ to $\Bbb R$, known as the dual space
of $X$) we write $H(X)$ instead of $H(w)$.
\end{defi}
\begin{defi}\label{def:rel} \cite[Definition 1.1.1]{DIM04} Let $X$ be
a subset of a Hausdorff topological group $(G,\tau).$  We say that
$X$ is  {\it relatively minimal} in $G$ if every coarser Hausdorff
group topology $\sigma \subset \tau$ of $G$ induces on  $X$ the
original topology. That is, $\sigma \big |_{X}=\tau \big |_{X}.$
\end{defi}
\begin{thm}\cite[Theorem 2.2]{MEG04}
\label{thm:rel} The subgroups  $X$ and $X^{\ast}$ are relatively
minimal in the generalized Heisenberg group $H(X)=(\Bbb R\times
X)\leftthreetimes X^{\ast}$ for every normed space $X$.
\end{thm}

The concept of co-minimality which is presented below played a major
role in generalizing and strengthen Theorem ~\ref{thm:rel}. Let $H$
be a subgroup of a topological group $(G,\gamma)$. The quotient
topology  on the left coset space $G/H:=\{gH\}_{g\in G}$ will be
denoted by $\gamma/H. $

\begin{defi}\label{def:com} \cite[Definition 1.1.2]{DIM04}
Let $X$ be a topological subgroup of a Hausdorff topological group
$(G,\tau)$ . We say that $X$ is {\it co-minimal} in $G$ if every
coarser Hausdorff group topology
 $\sigma \subset\tau$ of $G$ induces on the coset space $G/X$ the original topology.
That is, $\sigma/X=\tau/X$.
\end{defi}
\begin{defi} \label{d:min-maps}
Let $E, F,A$ be abelian Hausdorff groups. A biadditive mapping
$w:E\times F \rightarrow A$ will be called {\it separated} if for
every pair $(x_{0},f_0)$ of nonzero elements there exists a pair
$(x,f)$ such that $f(x_0)\neq 0_{A}$ and  $f_0(x)\neq 0_{A}$, where
$f(x)=w(x,f)$.\\
\end{defi}

\begin{defi} \label{def:stm}  \cite[Definition 2.2]{DIM04}
Let $(E,\sigma), (F,\tau),(A,\nu)$ be abelian  Hausdorff 
topological groups. A continuous separated biadditive mapping
$$w:(E,\sigma)\times (F,\tau) \rightarrow (A,\nu)$$ will be called
{\it strongly minimal }if for every coarser triple
$(\sigma_1,\tau_1,\nu_1)$ of Hausdorff group topologies
$\sigma_1\subset \sigma,\tau_1\subset\tau, \nu_1\subset\nu$ such
that
$$w:(E,\sigma_1)\times (F,\tau_1) \rightarrow (A,\nu_1)$$ is
continuous (in such cases we say that the triple
$(\sigma_1,\tau_1,\nu_1)$ is compatible) it follows that
$\sigma_1=\sigma,\tau_1=\tau$. We say that the biadditive mapping is
minimal if $\sigma_1=\sigma,\tau_1=\tau$ holds for every compatible
triple $(\sigma_1,\tau_1,\nu)$ (with $\nu_1:=\nu$).
\end{defi}
\begin{remark}\label{rem:stm}
The multiplication map $A \times A \to A$ is minimal for every
Hausdorff topological unital ring $A$. However note that the
multiplication map $\Z \times \Z \to \Z$ (being minimal) is not
strongly minimal.
\end{remark}

\vskip  0.5cm The following theorem which uses the concept of
co-minimality and strongly biadditive mappings generalizes Theorem
~\ref{thm:rel}.
\begin{thm} \label{thm:gen} \cite[Theorem 4.1]{DIM04}
Let $w:(E,\sigma)\times (F,\tau) \rightarrow (A,\nu)$ be a
strongly minimal biadditive mapping. Then: \\
1. $A,\ A\times E$ and $A\times F$ are co-minimal subgroups of the
Heisenberg group
$H(w)$.\\
2. $E\times F$ is a relatively minimal subset in $H(w)$. \\
3. The subgroups $E$ and $F$ are relatively minimal in $H(w)$.
\end{thm}
\begin{remark}
The mapping $w:X\times X^{\ast}\rightarrow \Bbb R$ $(x,f)\mapsto
f(x)$ is strongly minimal for every normed space $X$. Therefore,
Theorem ~\ref{thm:gen} is indeed a generalization of Theorem
~\ref{thm:rel}.
\end{remark}
\begin{corol}\cite[Corollary 4.2]{DIM04} \label{cor:ami}
The following conditions are equivalent:
\begin{enumerate}
\item $H(w)$ is a minimal group.
\item
$w$ is a minimal biadditive mapping and $A$ is a minimal group.
\end{enumerate}
\end{corol}
Since $\Bbb Z$ with the $p$-adic topology $\tau_p$  is a minimal
group
 for every prime $p$ \cite{PRO71} the following corollary is obtained by Remark
 ~\ref{rem:stm}:
\begin{corol} \cite[Corollary 4.6.2]{DIM04}
The Heisenberg group $H(w)=(\Z \times \Z) \leftthreetimes \Z$ of the
mapping $(\Z, \tau_p) \times (\Z, \tau_p) \to (\Z, \tau_p)$ is a
minimal two step nilpotent precompact group for every $p$-adic
topology $\tau_p$.
\end{corol}
 \vskip  0.5cm
\section{Topological rings and absolute values} \label{sec:abs}
\vskip 0.3cm
 In this paper  rings are not assumed to be
 necessarily
associative. However, when we consider division rings we assume they
are associative unless otherwise is stated.

\begin{defi} \label{rem:arc}
An absolute value $A$ on a (not necessarily associative) division
ring $K$ is {\it archimedean}  if there exists $n\in \Bbb N$ such
that
$A(n)>1$ (where, for any
$n\in \Bbb N, \ \
 n:=n.1=1+\cdots +1$ \ (n terms).
\end{defi}

\vskip 0.3cm From now on we use the following notations for a
commutative group $G$ which is denoted additively:
the zero element is denoted
  by $0_G$. If $G$ is also a ring with multiplicative unit we denote this element   by
  $1_G$.  In the case of a group $G$ which is a direct product
  of groups we shall use slightly different notation
 and denote the zero element by $\bar{0}_G$.

\begin{lemma} \label{l:unbounded nbd}
Let $X$ be a  (not necessarily associative) division ring with an
archimedean absolute value
 $A$ and denote by $\tau$ the ring topology induced by the
absolute value. Let $\sigma\subset \tau$ be a strictly coarser
 group topology with respect to the additive structure of
$X$. Then, every $\sigma$-neighborhood of $0_X$ is unbounded with
respect to the absolute value.
\end{lemma}
\begin{proof}
Since $\sigma$ is strictly coarser than $\tau$, there exists  an
open ball $B(0,r)$ with $r>0$ not containing any
$\sigma$-neighborhood of $0_X$. Then, for every
$\sigma$-neighborhood $U$ of $0_X$ there exists $x$ in $U$ such that
$A(x) \geq r$. Fix a $\sigma$-neighborhood $V$ of $0_X$. We show
that $V$ is unbounded with respect to the absolute value $A$. Since
$A$ is an archimedean absolute value there exists $n_0\in \Bbb N$
such that $A(n_0)=c>1.$ Clearly, for every $m\in \Bbb N$ there
exists a $\sigma$-neighborhood $W$ of $0_X$ such that
$$\underbrace{W+W+\cdots +W}_{n_{0}^{m}}\subset V.$$
By our assumption there exists $x\in W$ such that $A(x)\geq r$. Now
for the element
$$n_{0}^{m}x:= \underbrace{x+x+\cdots +x}_{n_{0}^{m}}\in V$$ we obtain
that $A(n_0^{m}x)=A(n_0)^{m} A(x) \geq c^{m}r$. This clearly
means that $V$ is unbounded.
\end{proof}

\begin{lemma} \label{l:prod}
Let $(G_i)_{i\in I}$ be a family of topological groups.  For each
$i\in I$ denote by $\tau_{i}$ the topology of $G_i$ and by $p_i$ the
projection of $G:=\prod_{i\in I}G_i$ to $G_i$. Suppose that $\sigma$
is a group topology on $G$ which is strictly coarser than the
product topology on $G$ denoted by $\tau$. Then there exist $j\in I$
and a group topology $\sigma_j$ on $G_j$ which is strictly coarser
than $\tau_j$, such that $\frak{B}_j=p_j(\frak{B})$, where
$\frak{B_j}$ is the neighborhood filter of $0_{G_j}$ with respect to
$\sigma_j$ and $\frak{B}$ is the neighborhood filter of $\bar{0}_G$
with respect to $\sigma$.
\end{lemma}
\begin{proof}
Since  the topology $\sigma$ is strictly coarser than $\tau$ which
is the product topology on $G$, we get that there exists $j\in I$
for which the projection $p_j:(G,\sigma)\rightarrow (G_j,\tau_j)$ is
not continuous at $\bar{0}_G$. Hence, there exist a
$\tau_j$-neighborhood $V$ of $0_{G_j}$ such that $p_j(O)\nsubseteq
V$ for every $O\in \frak{B}$. Hence, if $p_j(\frak{B})$ is the
neighborhood filter of $0_{G_j}$ for some group topology $\sigma_j$
on $G_j$ then this topology is strictly coarser than $\tau_j$. We
shall prove that this formulation defines a group topology
$\sigma_j$. Indeed, consider the normal subgroup $H=\prod_{i\in I}
F_i$ of $G$ where

$$F_i=
\left\{
  \begin{array}{ll}
    G_i & \textrm{if} \ \ i\neq j \\
  \{0_{G_{i}}\} & \textrm{if} \ \ i=j
  \end{array}
\right..$$

It is easy to show that $(G_j,\tau_j)$ is topologically isomorphic
to the quotient group $G/H$ of $(G,\tau)$.  Let $\sigma_j$ be the
finest topology on $G_j$ for which
 the projection $p_j:(G,\sigma)\rightarrow G_j$ is continuous. It is exactly the quotient
 topology on $G_j=G/H$ for the topological group $(G,\sigma)$. By
 our construction $\sigma_j$ is strictly coarser than $\tau_j$.
Then indeed $\sigma_j$ is the desired group topology on $G_j$ and
  $\frak{B}_j=p_j(\frak{B})$ is the desired neighborhood filter.
 \end{proof}
\begin{thm}
\label{thm:new} Let $F$ be a (not necessarily associative) division
ring furnished with an archimedean absolute value $A$. For each
$n\in \Bbb N$,
$$w_n:F^{n}\times F^{n}\mapsto F,\quad w_n(\bar{x},\bar{y})=\sum_{1=1}^{n}x_iy_i$$
(where $(\bar{x},\bar{y})=((x_1,\ldots x_n),(y_1,\ldots y_n))$ is a
strongly minimal biadditive mapping.
\end{thm}
\begin{proof}
Clearly, for each $n\in \Bbb N,\ w_n$ is a continuous separated
biadditive mapping.
 Denote by $\tau$ the topology of $F$ induced by $A$ and by
$\tau^{n}$ the product topology on $F^{n}$. Consider the max-metric
$d$ on $F^n$. Then its topology is exactly $\tau^{n}$.
 Let $(\sigma,\sigma',\nu)$ be a compatible triple with respect to
 $w_n$. We prove that $\sigma=\sigma'=\tau^{n}$. Assuming the
 contrary we get that at least one of the group topologies
 $\sigma,\sigma'$ is strictly coarser than $\tau^{n}$.
 We first assume that $\sigma$ is strictly coarser than $\tau^{n}.$
  Since $\nu$ is Hausdorff and $(\sigma,\sigma',\nu)$ is compatible there exist
a $\nu$-neighborhood $Y$ of $0:=0_F$ and $V,W$ which are
respectively $\sigma,\sigma'$-neighborhoods of  $\bar{0}_{F^{n}}$
such that $VW\subset Y$ and in addition $1_F\notin Y$. Since $W \in
\sigma'\subset \tau^{n}$, then there exists $\epsilon_0 >0$ such
that the corresponding $d$-ball $B(0,\epsilon_0)$ is a subset of
$W$. Since $\sigma$ is strictly coarser than $\tau^{n}$
(by Lemmas \ref{l:unbounded nbd} and \ref{l:prod})
 there exists $i\in I:=\{1,2, \cdots,n\}$ such that $p_i(V)$ is norm
 unbounded. Therefore, there exists $\bar{x}\in V$ such that
 $A(p_i(\bar{x}))>\frac{1}{\epsilon_0}$. Hence, $A((p_i(\bar{x}))^{-1})<\epsilon_0.$
  Now, let us consider a vector $\bar{a}\in  F^{n}$ such that for every $j\neq
 i,\ a_j=0$ and 
$a_i=(p_i(\bar{x}))^{-1}$. Clearly, $\bar{a}\in
B(0,\epsilon_0)\subset
 W$. We then get that $w_n(\bar{x},\bar{a})=1_F\in VW\subset Y$. This
 contradicts our assumption. Using the same technique we can show
 that $\sigma'$  can't be strictly coarser than $\tau^{n}.$
 \end{proof}
 \begin{example}
 \begin{enumerate}
\item Let $F\in \{\Bbb Q,\Bbb R,\Bbb C\}$ with the usual absolute
value. Then for each $n\in \Bbb N$ the map
$$w_n:F^{n}\times F^{n}\mapsto F$$ is strongly minimal. The case of $F$
equals to $\Bbb R$ follows also from
 \cite[Proposition 2.42]{DIM04}.

\item  For each $n\in \Bbb N$ the map $$w_n:\Bbb H^{n}\times \Bbb H^{n}\mapsto \Bbb H$$ is strongly minimal, where $\Bbb H$ is the quaternions ring equipped with the archimedean absolute
value defined by:
$$\|q\|=(a^{2}+b^{2}+c^{2}+d^{2})^{\frac{1}{2}}$$ for each
$q=a+bi+cj+dk\in \Bbb H.$

\item  Let $G$ be the non-associative ring of octonions.

This ring
 can be defined (see \cite{WIK07}) as pairs of quaternions (this is
the Cayley-Dickson construction). Addition is defined pairwise. The
product of two pairs of quaternions $(a, b)$ and $(c, d)$ is defined
by $(a, b)(c, d) = (ac - db^{\ast}, a^{\ast}d + cb)$ where
$z^{\ast}=e-fi-gj-hk$ denotes the conjugate of $z=e+fi+gj+hk.$ We
define a norm on $G$ as follows:$$
\|(a+bi+cj+dk,e+fi+gj+hk)\|=(a^{2}+b^{2}+c^{2}+d^{2}+e^{2}+f^{2}+g^{2}+h^{2})^{\frac{1}{2}}.$$
This norm agrees with the standard Euclidean norm on $\Bbb R^{8}$.
It can be proved that for each $x_1,x_2\in G, \
 \|x_1x_2\|=\|x_1\|\cdot \|x_2\|$ hence $\|\quad \|$ is an absolute
value and clearly it is archimedean. Again by Theorem \ref{thm:new}
the map $$w_n: G^{n}\times G^{n}\mapsto G$$ is strongly minimal for
each $n\in\Bbb N$.
\end{enumerate}
\end{example}
\begin{corol} \label{cor:nrs}
Under the conditions of Theorem ~\ref{thm:new}  we obtain the
following results:
\begin{enumerate}
\item $(F\times \{\bar{0}_{F^{n}}\})\leftthreetimes  \{\bar{0}_{F^{n}}\}, (F\times F^{n})\leftthreetimes
\{\bar{0}_{F^{n}}\}$ and $ (F\times
\{\bar{0}_{F^{n}}\})\leftthreetimes F^{n}$ are co-minimal subgroups
of the Heisenberg group $H(w_n)$.
\item $(\{0_{F}\}\times F^{n})\leftthreetimes F^{n}$ is a relatively
minimal subset in $H(w_n)$ .
\item The subgroups $(\{0_F\}\times F^{n})\leftthreetimes \{\bar{0}_{F^{n}}\}$ and $(\{0_F\}\times \{\bar{0}_{F^{n}}\})\leftthreetimes
F^{n}$ are relatively minimal in $H(w_n)$.
\end{enumerate}
\end{corol}
\begin{proof}
Apply Theorem  ~\ref{thm:gen} to the strongly minimal biadditive
mapping $w_n$.
\end{proof}
\begin{remark} \label{rem:nrs}
We replace $H(w_n)$ by $H(F^{n})$ for convenience ($w_n$ is the
strongly minimal biadditive mapping from ~\ref{thm:new}).
 In terms of matrices: $H(F^{n})$ is the
$2n+1$-dimensional Heisenberg group with coefficients from $F$ which
consists of square matrices of size $n+2$:
 $$A= \left(\begin{array}{ccccccc} 1_F & x_1 & x_2 & \ldots & x_{n-1} & x_{n} & r\\
0_F & 1_F & 0_F & 0_F & 0_F & 0_F & y_1\\
0_F & 0_F & \ddots & \ddots & \ddots & \vdots & y_2\\
\vdots & \vdots & \ddots & \ddots & \ddots & 0_F & \vdots\\
\vdots & \vdots & \ddots & \ddots & 1_F & 0_F & y_{n-1}\\
0_F & 0_F & \ddots & \ddots & 0_F & 1_F & y_n\\
0_F & 0_F & 0_F & \ldots & \ldots & 0_F & 1_F\\
 \end{array}\right)$$ and by the result (2) of Corollary ~\ref{cor:nrs}
 we obtain that the set of matrices  $$B= \left(\begin{array}{ccccccc} 1_F & x_1 & x_2 & \ldots & x_{n-1} & x_{n} & 0_F\\
0_F & 1_F & 0_F & 0_F & 0_F & 0_F & y_1\\
0_F & 0_F & \ddots & \ddots & \ddots & \vdots & y_2\\
\vdots & \vdots & \ddots & \ddots & \ddots & 0_F & \vdots\\
\vdots & \vdots & \ddots & \ddots & 1_F & 0_F & y_{n-1}\\
0_F & 0_F & \ddots & \ddots & 0_F & 1_F & y_n\\
0_F & 0_F & 0_F & \ldots & \ldots & 0_F & 1_F\\
 \end{array}\right)$$
 is a relatively minimal subset of $H(F^{n}).$
 \end{remark}
\begin{lemma} \label{lem:rms}
\begin{enumerate}
\item If $H$ is a subgroup of a topological group $(G,\tau)$ and $X$ is a
 relatively minimal subset in $H,$ then $X$ is also relatively minimal in  $G.$
\item  Let $(G_1,\tau_1),(G_2,\tau_2)$ be topological groups and $H_1,H_2$ be their
subgroups (respectively). If $H_1$ is relatively minimal in $G_1$
and there exists a topological isomorphism $f:
(G_1,\tau_1)\rightarrow (G_2,\tau_2)$ such that the restriction to
$H_1$ is a topological isomorphism onto $H_2,$ then $H_2$ is
relatively minimal in   $G_2.$
\item Let $(G,\tau)$ be a topological group and let $X$ be a
 subset of $G$. If $X$ is relatively minimal in $(G,\tau),$
 then every subset of $X$ is also relatively minimal in $(G,\tau)$.
\end{enumerate}
\end{lemma}
\begin{proof}
(1): Let $\sigma\subset \tau$ be a  coarser Hausdorff group topology
of $G,$ then  $\sigma\big |_H\subset\tau\big |_H$ is a coarser
Hausdorff group topology of $H.$ Since $X$ is a
 relatively minimal subset in $H,$ we get that $$\sigma\big |_X= (\sigma\big |_H)\big
 |_X=(\tau\big |_H)\big
 |_X=\tau\big |_X.$$ Hence, $X$ is relatively minimal in  $G.$\\
 (2): Observe that if $\sigma_2\subset \tau_2$ is a  coarser  Hausdorff group topology of $G_2,$
 then $$f^{-1}(\sigma_2)=\{f^{-1}(U)|\ U\in \sigma_2\}\subset
 \tau_1$$
is a  coarser group topology of $G_1.$ Since $H_1$ is relatively
minimal in $(G_1,\tau_1)$ we obtain that $\tau_1\big
|_{H_1}=f^{-1}(\sigma_2)\big |_{H_1}.$ This implies that $\tau_2\big
|_{H_2}=\sigma_2\big |_{H_2}.$ This completes our proof.\\
(3): Let $Y$ be a subset of $X$ and $\sigma\subset \tau$ a coarser
Hausdorff group topology. Then, by the fact that $X$ is relatively
minimal in $(G,\tau)$ and since $Y$ is a subset of $X$ we obtain
that $$\sigma\big |_Y=(\sigma\big |_X)\big
 |_Y=(\tau\big |_X)\big
 |_Y=\tau\big |_Y.$$ Hence, $Y$ is relatively minimal in $G.$
\end{proof}

The following is new even for the case of $F=\Bbb R$\ (for $n>1$).
\begin{thm}\label{thm:ne2}
Let  $F$ be a field furnished with an archimedean absolute value
$A$. For all $n\in \Bbb N$ denote by $U_{n+2}(F)$ the topological
group of all  $n+2\times n+2$ upper 
unitriangular matrices with entries from $F.$  Then $\forall n\in
\Bbb N$ and for each $i,j$ such that $i<j, \ (i,j)\neq (1,n+2)$ each
of the subgroups
 $$G^{n+2}_{ij}(F):=\Bigg \{\left(\begin{array}{ccccc} 1_F & 0_F & 0_F & 0_F & 0_F\\
0_F & 1_F & 0_F & 0_F & 0_F\\
\vdots & \ddots & \ddots & a_{ij} & \vdots\\
0_F & 0_F & 0_F & 1_F & 0_F\\
0_F & \ldots & \ldots & 0_F & 1_F\\
 \end{array}\right) \in U_{n+2}(F)\Bigg \}$$
 (where $a_{ij}$ is in the $ij$ entry)
 is relatively minimal in $U_{n+2}(F).$
\end{thm}
\begin{proof}
We   prove the assertion for two cases:
 First case: $i=1$ or $j=n+2$ (that is the indexes from the first row or from the last column)
  and the second case:  $i>1,\ j<n+2$.
 Let us consider the first case:
  we know by Remark  ~\ref{rem:nrs} that the set $S$ of square matrices of size $n+2$:
$$B= \left(\begin{array}{ccccccc} 1_F & x_1 & x_2 & \ldots & x_{n-1} & x_{n} & 0_F\\
0_F & 1_F & 0_F & 0_F & 0_F & 0_F & y_1\\
0_F & 0_F & \ddots & \ddots & \ddots & \vdots & y_2\\
\vdots & \vdots & \ddots & \ddots & \ddots & 0_F & \vdots\\
\vdots & \vdots & \ddots & \ddots & 1_F & 0_F & y_{n-1}\\
0_F & 0_F & \ddots & \ddots & 0_F & 1_F & y_n\\
0_F & 0_F & 0_F & \ldots & \ldots & 0_F & 1_F\\
 \end{array}\right)$$ is relatively minimal in $H(F^{n}).$
Since  $H(F^{n})$   is a subgroup of $U_{n+2}(F)$ we get by Lemma
~\ref{lem:rms} that $S$   is relatively minimal in $U_{n+2}(F).$
Now, $G^{n+2}_{1j}(F)\subset S$ for every $1<j<n+2$ and
$G^{n+2}_{in+2}(F)\subset S$ for every $1<i<n+2.$
 By Lemma ~\ref{lem:rms} we obtain that
$G^{n+2}_{ij}(F)$ is relatively minimal
 in $U_{n+2}(F)$ for every pair of indexes $(i,j)$ such that $i=1$
 or $j=n+2$ 
 (in addition to the demands: $i<j$ and$(i,j)\neq(1,n+2)$).
\vskip 0.5cm Case $2$: $i>1,\ j<n+2$.  Fix $n\in N$ and a pair
$(i,j)$ such that $1<i<j<n+2$.  We shall show that $G^{n+2}_{ij}(F)$
is relatively minimal in $U_{n+2}(F).$ We define the following
subgroup of
 $U_{n+2}(F)$: $$\tilde{U}_{n+2}(F):=\{A\in U_{n+2}(F)|\  a_{kl}=0_F \textrm{ if } l\neq k<i\}$$
 (it means that  the first $i-1$ rows of every matrix  contain  only  $0_F$ at each entry (besides the diagonal)).
 Clearly, this group is isomorphic to the group
 $U_{(n+2-(i-1))}(F)=U_{n+3-i}(F).$ Indeed, for every matrix $A\in
 \tilde{U}_{n+2}(F)$ if we delete the first $i-1$ rows and the first
 $i-1$ columns we obtain a matrix which belongs to $U_{n+3-i}(F)$ and
 it also clear that this way we obtain an isomorphism. Denote
 this isomorphism by $f.$
Now, $G^{n+2}_{ij}(F)$ is a subgroup of $\tilde{U}_{n+2}(F)$ and
$f(G^{n+2}_{ij}(F))=G^{n+3-i}_{1j+1-i}(F).$ Since  $1<i<j<n+2$ we
obtain that $i\leq n$ and hence $n+3-i\geq 3.$ Therefore, we can use
the reduction to case (1) to obtain that 
 $G^{n+3-i}_{1j+1-i}(F)$ is
relatively minimal in $U_{n+3-i}(F).$ By applying Lemma
~\ref{lem:rms} (with $G_1:=U_{n+3-i}(F), \ G_2:= \tilde{U}_{n+2}(F),
\ H_1:=G^{n+3-i}_{1j+1-i}(F)$ and $H_2:=G^{n+2}_{ij}(F)$) we can
conclude that $G^{n+2}_{ij}(F)$ is relatively minimal in
$\tilde{U}_{n+2}(F)$ and hence also in $U_{n+2}(F)$ which contains
$\tilde{U}_{n+2}(F)$ as a subgroup. This completes our proof.
\end{proof}
\begin{remark}
In the particular case of $F=\Bbb R$ we obtain by previous results
 that for every $n\in \Bbb N$ each of the subgroups
$G^{n+2}_{ij}(\Bbb R)$ is relatively minimal in $SL_{n+2}(\Bbb R $).
It is derived from the fact that  $SL_{m}(\Bbb R $) is minimal  for
every $m\in \Bbb N$ (see \cite{RES91, DIM04}). These groups are also
relatively minimal in $GL_{n+2}(\Bbb R $) which contains
$SL_{n+2}(\Bbb R $) as a subgroup (see Lemma ~\ref{lem:rms}).
Nevertheless, the fact that these groups are relatively minimal in
$U_{n+2}(\Bbb R $) is not derived  from the minimality of
$SL_{n+2}(\Bbb R $)  since $U_{n+2}(\Bbb R )$ is contained in
$SL_{n+2}(\Bbb R )$ and  not the opposite (that is $SL_{n+2}(\Bbb R
$) is not a subset of $U_{n+2}(\Bbb R )$).
\end{remark}
\begin{defi}
Let $K$ be a Hausdorff topological division ring. A topological
$K$-vector space $E$ is {\it straight } if $E$ is Hausdorff and for
every nonzero $c\in E, \lambda\rightarrow \lambda c$ is a
homeomorphism from $K$ to the one-dimensional subspace $Kc$ of $E$.
The Hausdorff topological division ring is {\it straight } if every
Hausdorff $K$-vector space is straight.
\end{defi}
\begin{thm}\cite[Theorem 13.8]{WAR93}
A nondiscrete locally retrobounded division ring is straight. In
particular, a division ring topologized by a proper absolute value
is straight.
\end{thm}
\begin{lemma}\label{lem:map}
Let $(F,\tau)$ be a unital Hausdorff topological ring.  Consider the
following cases: \begin{enumerate} \item $(F,\tau)$ is a minimal
topological group.
\item The multiplication map $w:(F,\tau)\times (F,\tau) \rightarrow
(F,\tau)$ is strongly minimal.
\item $(F,\tau)$ is minimal as a topological module over $(F,\tau)$
(i.e. there is  no strictly coarser Hausdorff topology $\sigma$ on
 $F$ for which $(F,\sigma)$ is a topological module over $(F,\tau)$).
 \item  $(F,\tau)$ is minimal as a topological ring (i.e. there is
 no strictly coarser Hausdorff ring topology on $F$).
 \end{enumerate}
Then: $$(1)\Rightarrow (2)\Rightarrow (3)\Rightarrow
 (4).$$
\end{lemma}
\begin{proof} $(1)\Rightarrow (2)$: If $F$ is a unital topological
ring then $w$ is
minimal. Indeed, let $(\sigma_1,\tau_1,\nu_1)$ be a compatible
triple then the identity maps $(F,\sigma_1)\rightarrow (F,\tau)$ and
$(F,\tau_1)\rightarrow (F,\tau)$ are continuous since  the
multiplication map $w:(F,\sigma_1)\times(F,\tau_1)\rightarrow
(F,\tau)$ is continuous at $(\lambda, 1_F), (1_F,\lambda)$ for every
$\lambda\in F$ and from the fact that $$\forall \lambda\in F \
w(\lambda, 1_F)=w(1_F,\lambda)=\lambda .$$ Clearly, in the case of a
minimal topological Hausdorff group the definition of a minimal
biadditive mapping and a strongly minimal biadditive mapping
coincide. The rest of the implications are trivial.
\end{proof}
\begin{remark}
Although $(1)\Rightarrow (2)$, the converse implication in general
is not true. For instance, the multiplication map $w:\Bbb R\times
\Bbb R \rightarrow \Bbb R $ is strongly minimal but $\Bbb R$ is not
minimal as a topological Hausdorff group.
\end{remark}
\begin{lemma}\label{lem:fin}
Let $(R,\tau)$ be a straight division ring. Let $\tau_0$ be a
strictly coarser Hausdorff topology on $\tau$. Then  $(R,\tau_0)$ is
not  a topological vector space over $(R,\tau)$.
\end{lemma}
\begin{proof}
Let $\tau_0\subset \tau$. We shall show that if $(R,\tau_0)$ is a
topological vector space then $\tau_0=\tau$.  In the definition of
straight division ring let $K=(R,\tau)$ and $E=(R,\tau_0)$ also let
$c=1$. Then it is clear that the identity mapping
$(R,\tau)\rightarrow (R,\tau_0)$ is a homeomorphism. Hence,
$\tau=\tau_0$.
\end{proof}

\begin{remark}
By our new results we get that in the case of archimedean absolute
value, conditions (2)-(4) of Lemma  ~\ref{lem:map}  hold. Since a
proper non-archimedean absolute valued division ring  is a straight
division ring we get by Lemma ~\ref{lem:fin} that    the conditions
(3)-(4) in Lemma  ~\ref{lem:map}  hold in this situation. The
question that remains open is whether the multiplication map
$$w:(F,\tau)\times (F,\tau)\rightarrow (F,\tau)$$ is strongly minimal
 where  $F$ is a division ring and the topology $\tau$ is induced by
a proper non-archimedean absolute value.
We ask even more concretely: is the multiplication map
$$w:\Bbb Q\times \Bbb Q\rightarrow \Bbb Q$$  strongly minimal when
$\Bbb Q$ is equipped with the $p$-adic topology?
\end{remark}
I would like to thank D. Dikranjan and M. Megrelishvili    for their
suggestions and remarks.
\bibliographystyle{amsalpha}

\end{document}